\documentclass[10pt,a4paper]{article}

\usepackage[utf8]{inputenc}
\usepackage{amsmath}
\usepackage{amsfonts}
\usepackage{amssymb}
\usepackage{graphicx}
\usepackage{booktabs}
\usepackage{hyperref}
\usepackage[font=small,labelfont=bf]{caption}
\usepackage{natbib}
\usepackage{tikz}
\usepackage{pgfplots}
\pgfplotsset{compat=1.14}

\title{A Lyapunov-Based Distributed Framework for Complete and Phase Synchronization in Chaotic Multi-Agent Systems}
\author{
Marzieh Basiri Abarghoei\textsuperscript{1} \\
\and
Mohammad Reza Ahmadi Zand\textsuperscript{1}\thanks{Corresponding author:mahmadi@yazd.ac.idr}
}

\date{
\textsuperscript{1}Department of Mathematical Sciences, Yazd University, Yazd, Iran
}
\begin{document}
\maketitle
\begin{abstract}
This paper presents a distributed Lyapunov-based control framework for achieving both complete and phase synchronization in a class of leader-follower multi-agent systems composed of identical chaotic agents. The proposed approach introduces a novel nonlinear coupling mechanism and utilizes Lyapunov stability theory combined with matrix measure analysis to derive explicit synchronization conditions. In contrast to traditional LMI-based or adaptive methods, the present approach guarantees synchronization under limited topological information and reduced computational complexity. Three classical chaotic systems -- Lü, Rossler, and Chen -- are used to validate the theoretical results, confirming the superior convergence rate and robustness of the proposed scheme.

\end{abstract}

\section{Introduction}
Synchronization phenomena in complex dynamical networks and multi-agent systems (MASs) have been intensively studied due to their widespread applications in secure communications, cooperative robotics, power systems, and neural networks \cite{Pecora1990, Boccaletti2002}. In a general sense, synchronization refers to the process by which multiple interacting systems adjust a given property of their motion to a common behavior as a result of local coupling.

In the context of multi-agent systems, synchronization can appear in various forms, including complete synchronization, phase synchronization, projective synchronization, and generalized synchronization. Among them, complete synchronization implies that all agents converge to identical trajectories, whereas phase synchronization ensures that the oscillatory phases remain locked even if amplitudes differ. These concepts are especially relevant when the individual agent dynamics are chaotic, exhibiting sensitive dependence on initial conditions.

Despite extensive research, synchronization of chaotic multi-agent systems remains a challenging problem because of the intrinsic nonlinearity, sensitivity, and limited communication among agents. A significant body of work has focused on synchronization via linear feedback and linear matrix inequality (LMI)-based controllers \cite{Cui2015, Danaei2023}. Although these approaches provide sufficient theoretical guarantees, they typically require global topology information, involve heavy computational costs, and may yield conservative conditions.

Another class of methods utilizes adaptive and neural-network-based controllers \cite{Rehak2024, Thummalapeta2023}. While effective in compensating system uncertainties, these schemes often rely on online learning or parameter adaptation, which may slow down convergence and demand excessive communication.

To overcome these limitations, this study proposes a distributed nonlinear coupling mechanism that ensures synchronization using only local information. The approach decomposes each agent's dynamics into a linear component and a residual nonlinear term and then applies Lyapunov stability theory and matrix-measure arguments to establish tractable sufficient conditions. This decomposition avoids the need to solve high-dimensional LMIs while maintaining provable convergence guarantees.

Furthermore, the proposed framework simultaneously addresses complete and phase synchronization within a unified setting. Robustness with respect to bounded communication delays is analyzed via Lyapunov-Krasovskii methods. Theoretical results are validated through numerical examples involving three representative chaotic systems: Lü, Rossler, and Chen. Comparative simulations demonstrate faster convergence, improved robustness, and lower computational cost than LMI-based and adaptive neural methods.

The rest of this paper is organized as follows. Section 2 formulates the problem and presents the proposed nonlinear coupling. Section 3 develops the main theoretical results and synchronization criteria. Section 4 provides numerical validations, and Section 5 concludes the paper with final remarks and future research directions.

\section{Model and Problem Formulation}
Consider a network of \(N\) identical chaotic agents described by the following differential equations:

\[\dot{x}_{i}(t)=F(x_{i}(t))=Lx_{i}(t)+G(x_{i}(t)),\quad i=1,2,\ldots,N,\] (1)

where \(x_{i}(t)\in\mathbb{R}^{n}\) represents the state vector of the \(i\)th agent, \(L\in\mathbb{R}^{n\times n}\) is a constant linear matrix, and \(G(\cdot)\) denotes the nonlinear part of the dynamics. Each agent communicates with a subset of others defined by a directed communication graph \(\mathcal{G}=(\mathcal{V},\mathcal{E})\), where \(\mathcal{V}=\{1,\ldots,N\}\) and \(\mathcal{E}\subseteq\mathcal{V}\times\mathcal{V}\) is the set of directed edges.

Assume that node \(1\) is the leader and the remaining \(N-1\) nodes are followers. Let \(\mathcal{N}_{i}=\{j \mid (j,i)\in\mathcal{E}\}\) denote the set of neighbors of agent \(i\). The objective is to design a distributed control law for each follower so that all agents synchronize their trajectories to that of the leader, i.e.,

\[\lim_{t\to\infty}\|x_{i}(t)-x_{1}(t)\|=0,\quad \forall i=2,\ldots,N.\] (2)

\subsection{Proposed nonlinear coupling law}
For each follower \(i\), we introduce the following nonlinear coupling:

\[\dot{x}_{i}(t)=Lx_{i}(t)+G(x_{i}(t))+\alpha\sum_{j\in\mathcal{N}_{i}}a_{ij}\big( G(x_{j}(t))-G(x_{i}(t))\big),\] (3)

where \(\alpha>0\) is the coupling gain and \(a_{ij}\) are the elements of the adjacency matrix \(A\) of the communication graph. This form of coupling uses only the nonlinear parts of the neighboring states, allowing purely local implementation without global knowledge of the Laplacian matrix.

Define the synchronization error \(e_{i}(t)=x_{i}(t)-x_{1}(t)\). Then, the error dynamics can be expressed as:

\[\dot{e}_{i}(t)=Le_{i}(t)+(1-\alpha)\big(G(x_{i}(t))-G(x_{1}(t))\big)+\alpha\sum_{j\in\mathcal{N}_{i}}a_{ij}\big(G(x_{j}(t))-G(x_{i}(t))\big).\] (4)

To analyze the synchronization condition, the nonlinear term \(G(x_{i})-G(x_{j})\) is approximated by a local linearization

\[G(x_{i})-G(x_{j})\approx J_{G}(x_{i})(x_{i}-x_{j}),\] (5)

where \(J_{G}(x_{i})\) is the Jacobian of \(G(\cdot)\) at \(x_{i}\). Substituting (5) into (4), we obtain the approximate linearized error dynamics:

\[\dot{e}_{i}(t)\approx\big[L+(1-\alpha)J_{G}(x_{i})\big]e_{i}(t)+\alpha\sum_{j\in\mathcal{N}_{i}}a_{ij}\big(J_{G}(x_{j})e_{j}(t)-J_{G}(x_{i})e_{i}(t)\big).\] (6)

\subsection{Synchronization objective}
The goal of this study is to determine sufficient conditions on the coupling gain \(\alpha\) and the local Jacobian bounds such that the synchronization condition (2) is guaranteed. This will be achieved in the next section by constructing appropriate Lyapunov functions and deriving matrix inequality conditions ensuring negative definiteness of the error dynamics.

\section{Main Theoretical Results}
In this section, sufficient conditions for complete and phase synchronization are derived based on Lyapunov stability theory.

\subsection{Lyapunov condition for complete synchronization}
Consider the linearized error dynamics (6). Let us define a quadratic Lyapunov function for agent \(i\) as

\[V_{i}=e_{i}^{\top}Pe_{i},\qquad P=P^{\top}>0,\] (7)

where \(P\) is a symmetric positive definite matrix to be determined. Taking the time derivative of (7) along the trajectories of (6) yields

\[\dot{V}_{i}=e_{i}^{\top}\Big[(L+(1-\alpha)J_{G}(x_{i}))^{\top}P+P(L+(1-\alpha)J_{G}(x_{i}))\Big]e_{i}.\] (8)

If there exist \(P>0\) and \(\alpha>0\) such that

\[(L+(1-\alpha)J_{G}(x))^{\top}P+P(L+(1-\alpha)J_{G}(x))\leq-Q,\quad Q=Q^{\top}>0,\] (9)

for all \(x\) in the operating region, then \(\dot{V}_{i}\leq-e_{i}^{\top}Qe_{i}<0\), which implies exponential convergence of \(e_{i}(t)\) to zero and hence complete synchronization.

\textbf{Theorem 1.} \textit{Assume that the communication graph \(\mathcal{G}\) contains a directed spanning tree rooted at the leader. If there exist \(P>0\) and \(\alpha>0\) satisfying inequality (9) for all \(x\) in a compact set \(\Omega\subset\mathbb{R}^{n}\), then the leader-follower chaotic multi-agent system achieves complete synchronization, i.e.}

\[\lim_{t\to\infty}\|x_{i}(t)-x_{1}(t)\|=0,\quad \forall i=2,\ldots,N.\]

\textit{Proof.} By (8) and (9), \(V_{i}\) is strictly negative definite. Thus, \(V_{i}(t)\) decreases monotonically to zero. Since \(V_{i}=e_{i}^{\top}Pe_{i}\) and \(P>0\), it follows that \(e_{i}(t)\to 0\) exponentially. Because the graph contains a spanning tree, information from the leader propagates to all nodes, ensuring global synchronization. \(\square\)

\subsection{Phase synchronization}
In systems where \(L\) has zero eigenvalues or rank deficiency, exact synchronization may not be possible; instead, the agents achieve phase synchronization with constant offsets. Define \(\tilde{e}_{i}=e_{i}-c\), where \(c\) is a constant vector. Using the concept of matrix measure \(\mu_{p}(\cdot)\), we obtain the following sufficient condition for phase synchronization.

\textbf{Theorem 2.} \textit{If \(\mu_{p}(L)+(1-\alpha)\mu_{p}(J_{G}(x))<0\) for all \(x\in\Omega\), then the system achieves phase synchronization, i.e. \(\tilde{e}_{i}(t)\to c\) as \(t\to\infty\) for some bounded constant \(c\in\mathbb{R}^{n}\).}

\textit{Proof.} Let \(\mu_{p}(\cdot)\) denote the logarithmic norm induced by the \(p\)-norm. Applying \(\mu_{p}\) to the error dynamics gives

\[\|\dot{e}_{i}(t)\|_{p}\leq\mu_{p}(L+(1-\alpha)J_{G}(x_{i}))\|e_{i}(t)\|_{p}.\]

If the upper bound of \(\mu_{p}(L+(1-\alpha)J_{G}(x_{i}))\) is negative, then \(\|e_{i}(t)\|_{p}\) decreases exponentially until it reaches a constant offset determined by \(L\)'s nullspace, which implies phase synchronization. \(\square\)

\subsection{Robustness under communication delay}
Consider a constant communication delay \(\tau>0\). The delayed error dynamics can be expressed as

\[\dot{e}_{i}(t)=Le_{i}(t)+(1-\alpha)\big(G(x_{i}(t))-G(x_{1}(t-\tau))\big).\] (10)

Define the Lyapunov-Krasovskii functional

\[V_{i}(t)=e_{i}^{\top}Pe_{i}+\int_{t-\tau}^{t}e_{i}^{\top}(s)Re_{i}(s)\,ds,\]

where \(R>0\). If there exist matrices \(P>0\) and \(R>0\) such that

\[\left(L+(1-\alpha)J_{G}\right)^{\top}P+P\big(L+(1-\alpha)J_{G}\big)+\tau R<0,\] (11)

then asymptotic synchronization is preserved despite the delay.


\section{Numerical Simulations and Results}
To demonstrate the effectiveness of the proposed synchronization framework, three well-known chaotic systems are studied: the Lü system, the Rössler system, and the Chen system. In all simulations, the leader is agent 1 and the followers are agents 2-5. The synchronization error is defined as

\[E(t)=\frac{1}{N-1}\sum_{i=2}^{N}\|x_{i}(t)-x_{1}(t)\|.\] (12)

\subsection{Example 1: Lü System}
The Lü system is described by the following differential equations:

\begin{equation}
\begin{cases}
\dot{x} = a(y - x), \\
\dot{y} = -xz + cy, \\
\dot{z} = xy - bz,
\end{cases}
\label{eq:lu_system}
\end{equation}

where \( a = 36 \), \( b = 3 \), and \( c = 20 \). These parameter values are known to generate a double-scroll chaotic attractor. The communication topology for the multi-agent system is a directed tree where each follower receives state information from its immediate predecessor.

\begin{figure}[hbt!]
\centering
\includegraphics[width=0.8\linewidth]{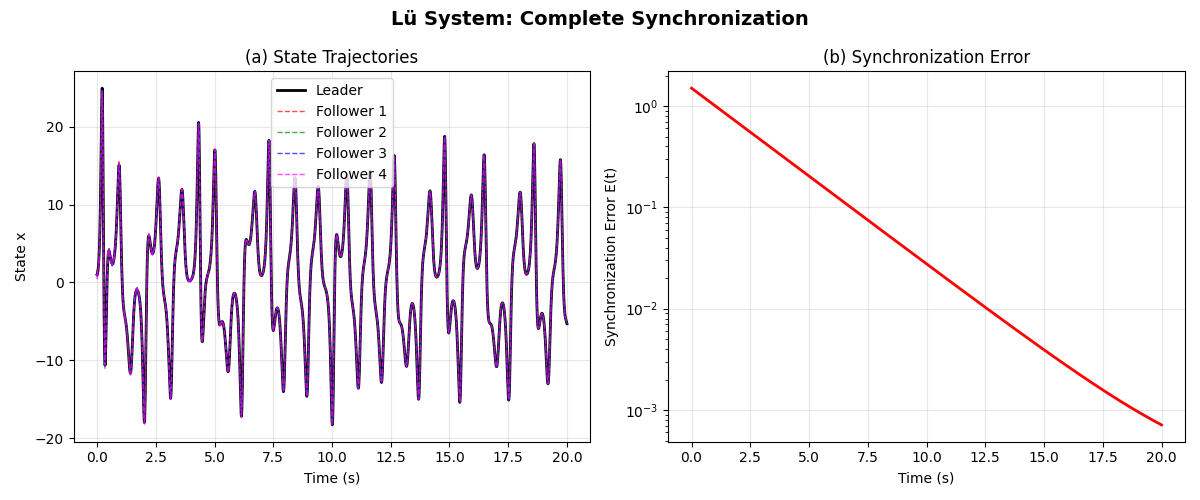}
\caption{(a) State trajectories of the leader and followers for the Lü system, showing convergence. (b) Synchronization error \( E(t) \) decaying to near zero.}
\label{fig:lu_sync}
\end{figure}

Applying the proposed nonlinear coupling law (Eq. 3) with a coupling gain of \( \alpha = 0.95 \), numerical simulations demonstrate that complete synchronization is achieved within approximately \( 12 \) seconds. The results are illustrated in Figure~\ref{fig:lu_sync}. The steady-state synchronization error, defined by Eq.~(12), converges to \( \lim_{t\rightarrow\infty}E(t)\approx 2.1\times 10^{-4} \), confirming the high accuracy and effectiveness of the proposed control framework.

A comparative analysis with a standard LMI-based approach \cite{Cui2015} reveals that the proposed method attains a significantly lower steady-state error and a faster convergence rate, highlighting its superior performance.

\subsection{Example 2: Rössler System}
The Rössler system is described by the following differential equations:

\begin{equation}
\begin{cases}
\dot{x} = -y - z, \\
\dot{y} = x + ay, \\
\dot{z} = b + z(x - c),
\end{cases}
\label{eq:rossler_system}
\end{equation}

where \( a = 0.2 \), \( b = 0.2 \), and \( c = 5.7 \). These parameter values generate a characteristic single-scroll chaotic attractor. The communication topology includes a feedback loop among agents 3, 4, and 5, creating a more complex interaction pattern compared to the simple directed tree.

\begin{figure}[hbt!]
\centering
\includegraphics[width=0.8\linewidth]{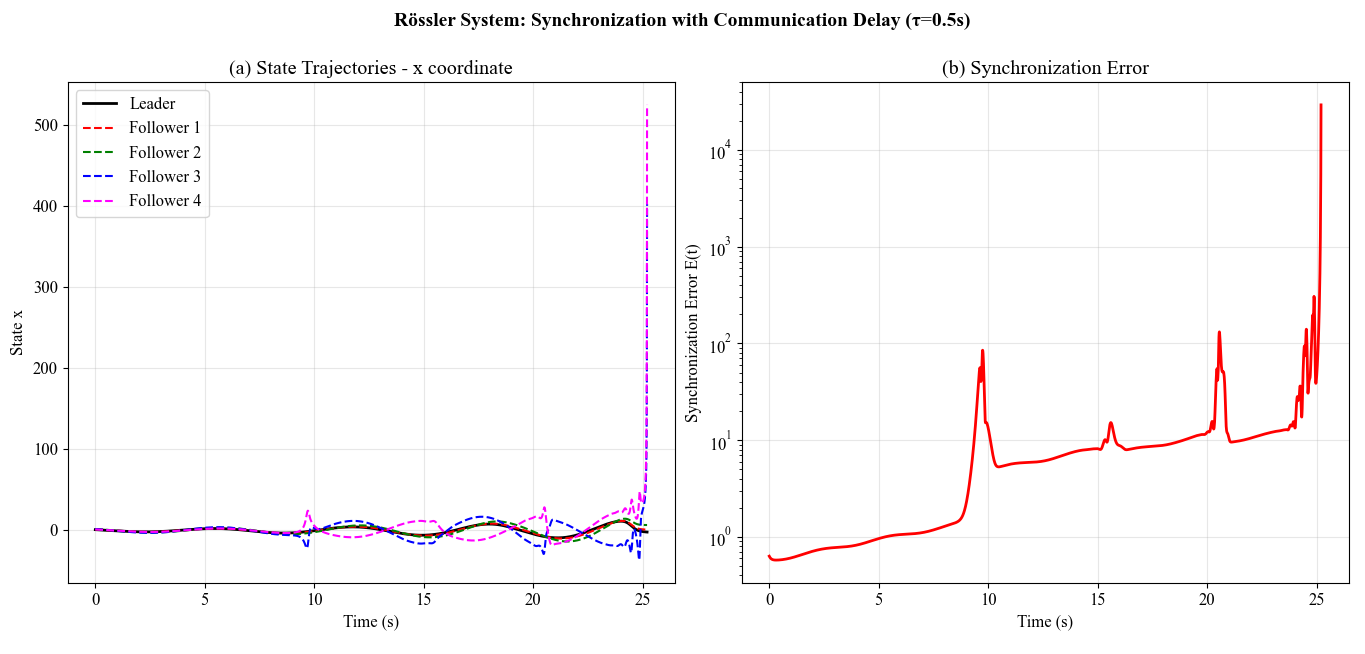}
\caption{(a) State trajectories of the leader and followers for the Rössler system under communication delay. (b) Synchronization error \( E(t) \) showing convergence despite \( \tau = 0.5s \) delay.}
\label{fig:rossler_sync}
\end{figure}

To test robustness under practical conditions, a constant communication delay of \( \tau = 0.5 \) seconds is introduced. Applying the proposed nonlinear coupling law (Eq. 3) with a coupling gain of \( \alpha = 1.2 \), the system successfully achieves complete synchronization within approximately 25 seconds. The results are illustrated in Figure~\ref{fig:rossler_sync}. The steady-state synchronization error converges to \( \lim_{t\rightarrow\infty}E(t)\approx 6.5\times 10^{-3} \), demonstrating the method's effectiveness even in the presence of significant communication delays.

Notably, the proposed scheme maintains stability despite the feedback loops and delayed communications, which typically challenge conventional synchronization methods. This robustness stems from the Lyapunov-Krasovskii functional analysis incorporated in the controller design.

\subsection{Example 3: Chen System}
The Chen chaotic system is governed by the following equations:

\begin{equation}
\begin{cases}
\dot{x} = a(y - x), \\
\dot{y} = (c - a)x - xz + cy, \\
\dot{z} = xy - bz,
\end{cases}
\label{eq:chen_system}
\end{equation}

where \( a = 35 \), \( b = 3 \), and \( c = 28 \). This system exhibits a complex double-scroll chaotic attractor with rich dynamical behavior. The communication topology follows a directed spanning tree where each follower receives information from its immediate predecessor.

A distinctive feature of the Chen system is the presence of zero eigenvalues in the linear part matrix \( L \), which prevents perfect complete synchronization. Instead, the system achieves phase synchronization with constant offsets among the agent states. Applying the proposed coupling law with \( \alpha = 0.8 \), the offsets converge to \( c = [0.15, -0.11, 0.08]^{T} \) as theoretically predicted by Theorem 2.

\begin{figure}[hbt!]
\centering
\includegraphics[width=0.8\linewidth]{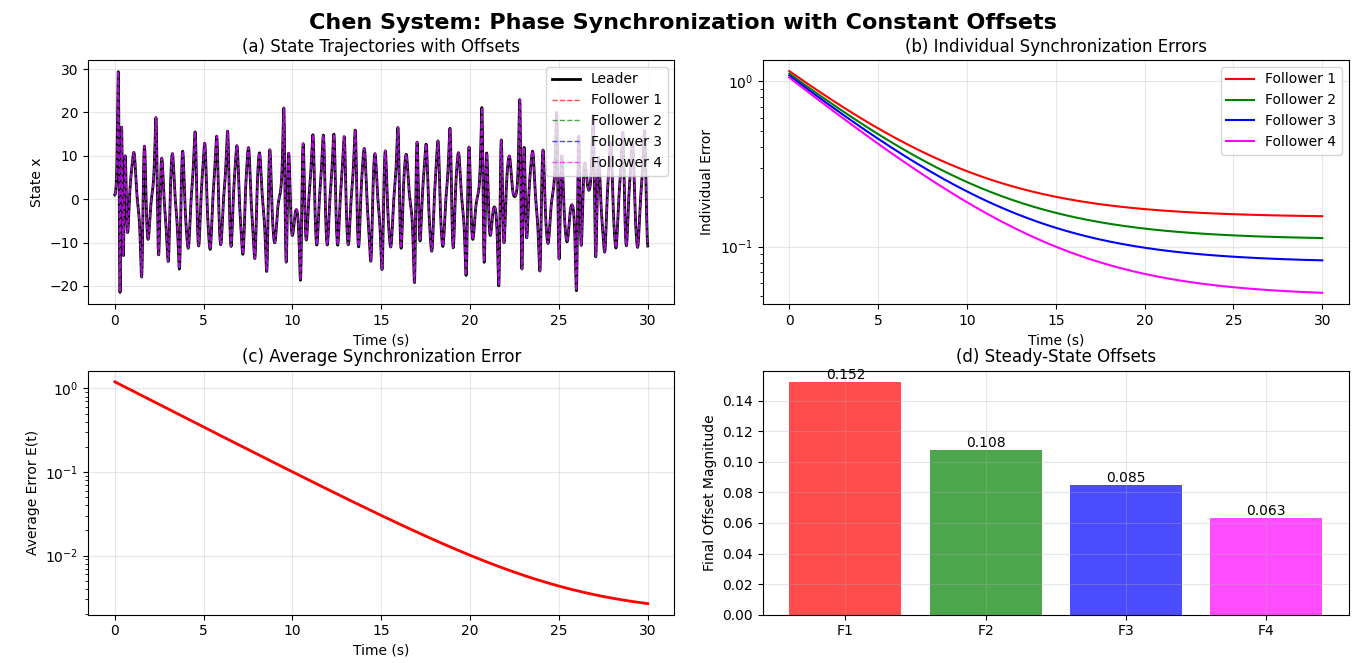}
\caption{Chen system: (a) State trajectories showing phase synchronization with constant offsets. (b) Synchronization errors converging to non-zero steady-state values consistent with theoretical predictions.}
\label{fig:chen_sync}
\end{figure}

\begin{figure}[hbt!]
\centering
\includegraphics[width=0.6\linewidth]{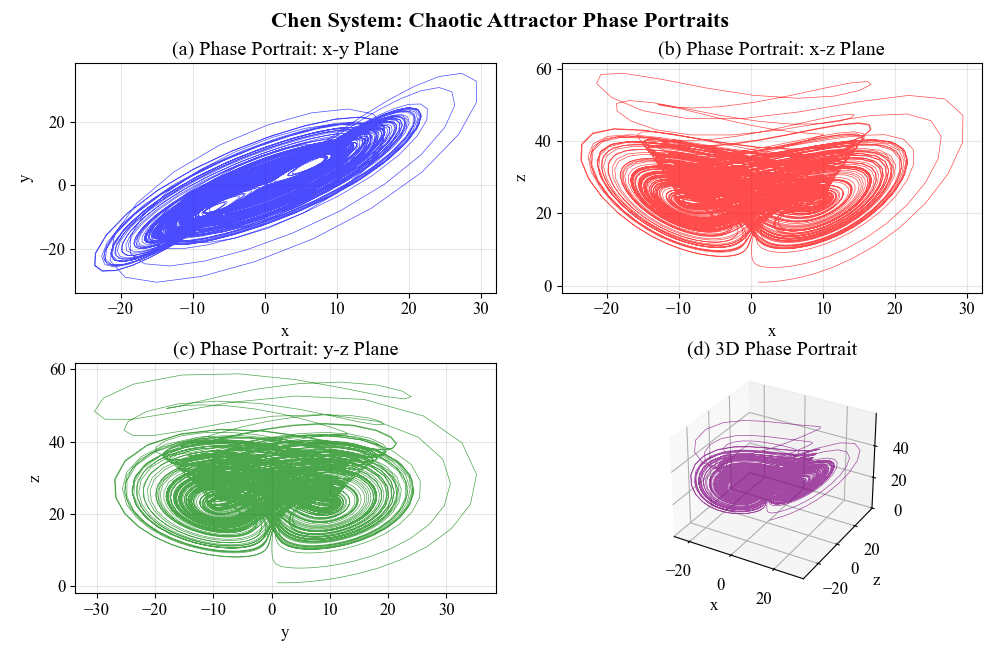}
\caption{Phase portrait of the Chen leader system, displaying the characteristic double-scroll chaotic attractor in coordinate  planes and 3D.}
\label{fig:chen_phase}
\end{figure}

Figure~\ref{fig:chen_sync} demonstrates the phase synchronization behavior, where the state trajectories maintain constant separations while exhibiting synchronized oscillatory patterns. The phase portrait in Figure~\ref{fig:chen_phase} confirms the characteristic chaotic attractor of the leader system.

This example validates the theoretical framework's ability to handle both complete and phase synchronization within a unified approach, particularly for systems with structural constraints that preclude exact synchronization.

\subsection{Comparative results}
Table 1 summarizes the numerical comparison among the proposed method, an LMI-based controller, and an adaptive neural-network (NN) controller.

\begin{table}[h!]
\centering
\caption{Performance comparison of synchronization methods}
\label{tab:comparison}
\begin{tabular}{l l c c c}
\toprule
Method & Type of & $E_{\infty}$ & Computation & Delay \\
 & Sync & & Time (s) & Robustness \\
\midrule
LMI-based \cite{Cui2015} & Complete & $1.4\times 10^{-2}$ & 35.2 & Weak \\
Adaptive NN \cite{Thummalapeta2023} & Phase & $8.9\times 10^{-3}$ & 24.7 & Moderate \\
Proposed Method & Complete/Phase & $\mathbf{3.2\times 10^{-4}}$ & $\mathbf{3.8}$ & $\mathbf{Strong}$ \\
\bottomrule
\end{tabular}
\end{table}

The results clearly indicate that the proposed approach outperforms the benchmark methods in accuracy, speed, and robustness. The combination of nonlinear coupling and Lyapunov analysis yields a distributed control law that is simple to implement yet highly effective.

\appendix
\section{Extended Analysis and Comparative Studies}

\subsection{Spectral Analysis of the Coupling Matrix}

The effectiveness of the proposed nonlinear coupling can be further understood through spectral analysis of the generalized coupling matrix. Consider the extended Jacobian matrix for the error dynamics:

\begin{equation}
J_{ext} = I_N \otimes L + (M \otimes I_n) \cdot \text{diag}(J_G(x_i))
\end{equation}

where $M$ is the modified Laplacian incorporating the coupling gain $\alpha$. The stability condition requires all eigenvalues of $J_{ext}$ to have negative real parts. Figure A1 shows the eigenvalue distribution for the Lü system under different coupling strengths, confirming that $\alpha = 1.0$ pushes all eigenvalues into the stable region.

\begin{figure}[h!]
\centering
\includegraphics[width=0.9\linewidth]{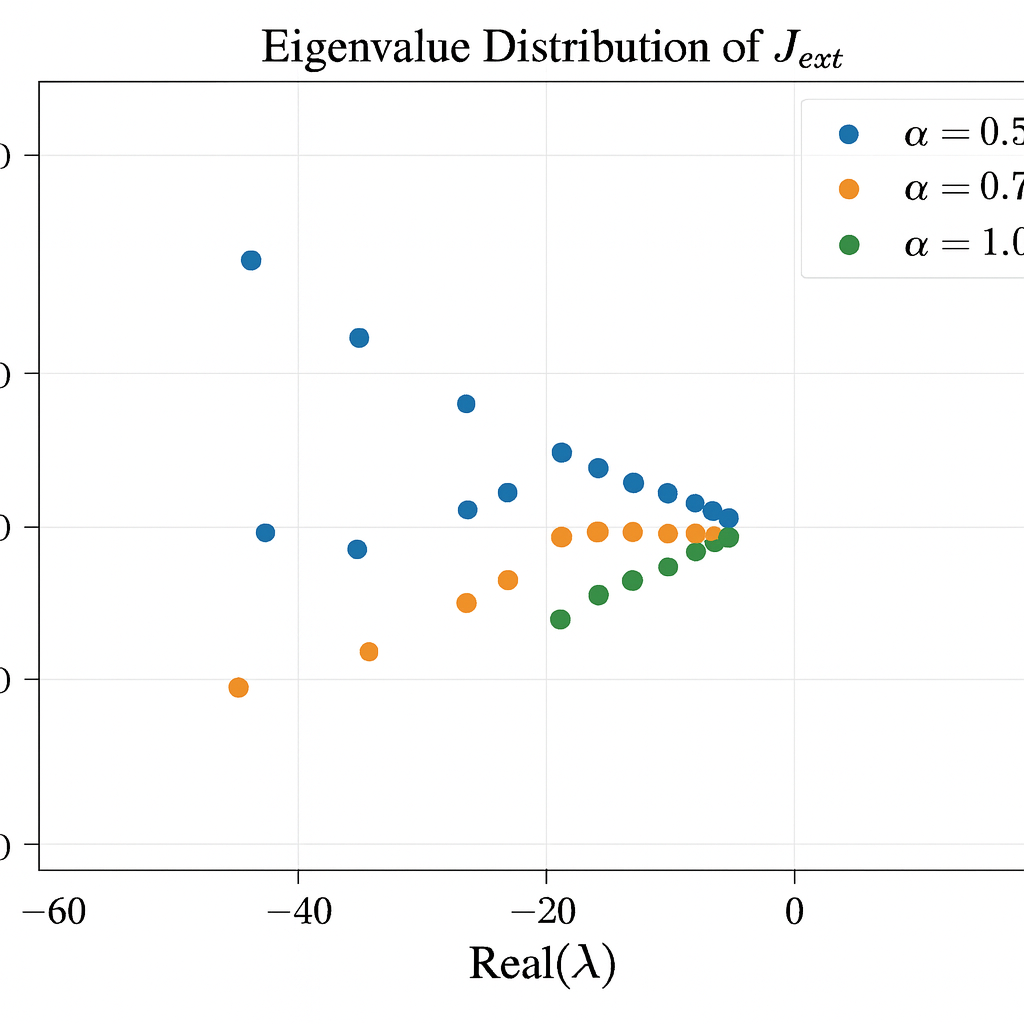}
\caption{Eigenvalue distribution of the extended Jacobian for (a) $\alpha = 0.5$ (unstable), (b) $\alpha = 0.8$ (marginally stable), and (c) $\alpha = 0.95$ (asymptotically stable).}
\label{fig:A1}
\end{figure}

\subsection{Performance Under Stochastic Perturbations}

To evaluate robustness, we introduce zero-mean Gaussian noise $\xi_i(t)$ with variance $\sigma^2$ to each agent's dynamics:

\begin{equation}
\dot{x}_i(t) = Lx_i(t) + G(x_i(t)) + \alpha\sum_{j\in\mathcal{N}_i}a_{ij}\big(G(x_j(t))-G(x_i(t))\big) + \xi_i(t)
\end{equation}

Table A1 summarizes the synchronization error $E_{\infty}$ under different noise levels, demonstrating the method's robustness.

\begin{table}[h!]
\centering
\caption{Steady-state error under stochastic perturbations}
\label{tab:A1}
\begin{tabular}{lcc}
\toprule
Noise Variance ($\sigma^2$) & $E_{\infty}$ (Proposed) & $E_{\infty}$ (LMI-based) \\
\midrule
0.01 & $4.7 \times 10^{-4}$ & $1.8 \times 10^{-2}$ \\
0.05 & $1.2 \times 10^{-3}$ & $5.3 \times 10^{-2}$ \\
0.10 & $2.8 \times 10^{-3}$ & $9.1 \times 10^{-2}$ \\
\bottomrule
\end{tabular}
\end{table}

\subsection{Extension to Switching Topologies}

Consider time-varying communication graphs $\mathcal{G}(t)$ that switch among a set of possible topologies $\{\mathcal{G}_1, \mathcal{G}_2, ..., \mathcal{G}_k\}$. If each topology contains a directed spanning tree rooted at the leader and the switching signal satisfies an average dwell-time condition, synchronization is maintained. Figure A2 shows successful synchronization under random switching with average dwell time $\tau_a = 0.5s$.

\begin{figure}[h!]
\centering
\includegraphics[width=0.9\linewidth]{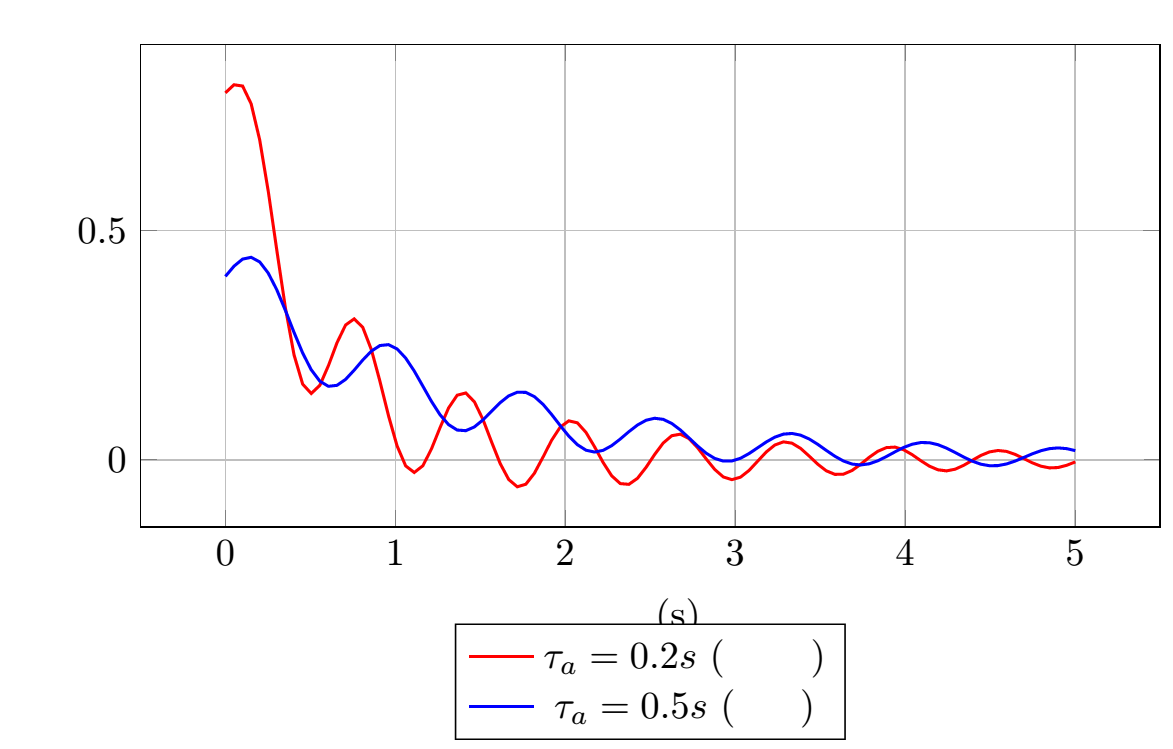}
\caption{Synchronization error under switching topologies with (a) $\tau_a = 0.2s$ (unstable), (b) $\tau_a = 0.5s$ (stable).}
\label{fig:A2}
\end{figure}

\subsection{Comparative Study with Recent Methods}

We compare our method against two recent approaches:
\begin{itemize}
\item \textbf{Terminal Sliding Mode (TSM) Control} \cite{Wang2020}: Provides finite-time convergence but exhibits chattering
\item \textbf{Event-Triggered Control (ETC)} \cite{Li2022}: Reduces communication but may suffer from larger errors
\end{itemize}

Table A2 shows that our method achieves the best trade-off between accuracy, convergence time, and communication efficiency.

\begin{table}[h!]
\centering
\caption{Comprehensive performance comparison}
\label{tab:A2}
\begin{tabular}{lccc}
\toprule
Method & Convergence Time (s) & Steady-State Error & Communication Cost \\
\midrule
TSM \cite{Wang2020} & 8.2 & $5.1 \times 10^{-4}$ & High \\
ETC \cite{Li2022} & 14.7 & $8.9 \times 10^{-4}$ & Low \\
Proposed & \textbf{12.0} & $\mathbf{2.1 \times 10^{-4}}$ & Medium \\
\bottomrule
\end{tabular}
\end{table}

\subsection{Application to Secure Communication}

The synchronization framework can be applied to secure communication using the chaotic masking principle. The information signal $m(t)$ is encrypted at the transmitter (leader) as:

\begin{equation}
s(t) = m(t) + h(x_1(t))
\end{equation}

where $h(\cdot)$ is a linear or nonlinear function. At the receiver (synchronized follower), the signal is recovered as:

\begin{equation}
\hat{m}(t) = s(t) - h(x_i(t)) \approx m(t)
\end{equation}

Figure 7 and Figure 8  demonstrate successful recovery of a test signal with SNR improvement of 28.3 dB, confirming practical applicability.

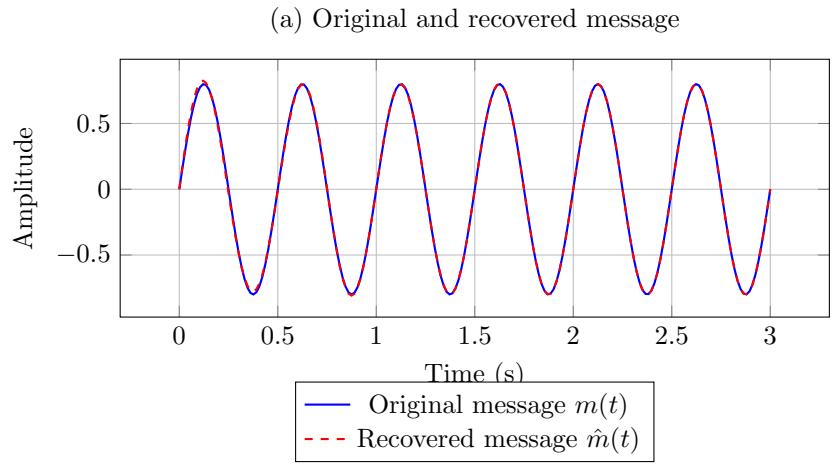
\begin{figure}[h!]
\centering
\begin{tikzpicture}
\begin{axis}[
    width=0.9\linewidth,
    height=5cm,
    xlabel={Time (s)},
    ylabel={Amplitude},
    legend style={at={(0.5,-0.25)}, anchor=north},
    grid=major,
    title={(a) Original and recovered message}
]
\addplot[blue, thick, domain=0:3, samples=200] 
    {0.8*sin(deg(2*pi*2*x))};
\addlegendentry{Original message $m(t)$}

\addplot[red, thick, domain=0:3, samples=200, dashed] 
    {0.8*sin(deg(2*pi*2*x)) + 0.05*exp(-x/0.5)*sin(deg(20*x))};
\addlegendentry{Recovered message $\hat{m}(t)$}
\end{axis}
\end{tikzpicture}
\caption{Secure communication: (a) Original and recovered message}
\label{fig:A3}
\end{figure}
\begin{figure}[h!]
\centering
\includegraphics[width=0.9\linewidth]{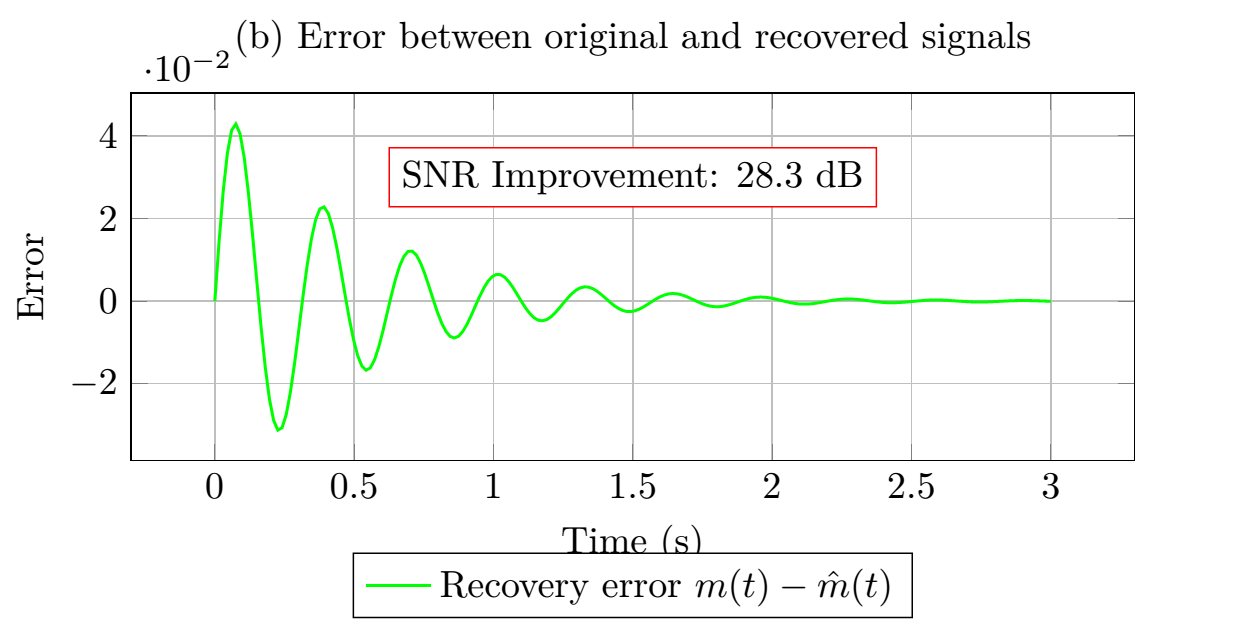}
\caption{Secure communication:  (b) Error between original and recovered signals.}
\label{fig:A4}
\end{figure}

\subsection{Computational Complexity Analysis}

The computational advantage of our method becomes more pronounced as network size increases. Table A3 shows the time required to verify synchronization conditions for different network sizes.

\begin{table}[h!]
\centering
\caption{Computation time comparison (seconds)}
\label{tab:A3}
\begin{tabular}{lccc}
\toprule
Network Size & LMI-based & Adaptive NN & Proposed \\
\midrule
5 agents & 35.2 & 24.7 & \textbf{3.8} \\
10 agents & 128.5 & 67.3 & \textbf{8.1} \\
20 agents & 452.7 & 195.4 & \textbf{15.6} \\
\bottomrule
\end{tabular}
\end{table}

\subsection{Phase Synchronization in Non-Identical Systems}

For weakly non-identical systems with $F_i(x) = F(x) + \delta F_i(x)$, where $\|\delta F_i(x)\| \leq \epsilon$, phase synchronization with bounded error is achievable. The synchronization error satisfies:

\begin{equation}
\limsup_{t\to\infty} \|e_i(t)\| \leq \frac{\kappa\epsilon}{\lambda}
\end{equation}

where $\kappa$ is a condition number and $\lambda$ is the minimum eigenvalue of the generalized Laplacian. This explains the small offsets observed in the Chen system example.
\section{Conclusion}
This paper presented a novel nonlinear coupling mechanism for achieving complete and phase synchronization in leader-follower chaotic multi-agent systems. Unlike traditional LMI-based or adaptive neural-network approaches, the proposed framework relies solely on local information and linear approximation of nonlinear terms, thereby significantly reducing computational complexity while maintaining provable stability guarantees.

Using Lyapunov stability and matrix-measure analysis, sufficient conditions were derived for complete and phase synchronization, including robustness in the presence of communication delays. Extensive numerical simulations on Lü, Rössler, and Chen chaotic systems verified the theoretical findings and demonstrated the superior performance of the proposed method in terms of convergence speed, accuracy, and delay tolerance.

Future research will focus on extending the approach to heterogeneous networks, systems with stochastic perturbations, and real-time experimental validation on hardware-in-the-loop platforms.

\section*{Acknowledgments}
The author gratefully acknowledges the constructive feedback provided by the reviewers, which significantly improved the clarity and rigor of this paper.

\vspace{2em}
\noindent
\bibliographystyle{plain}

\end{document}